\documentclass{amsart}
\usepackage{amsmath,amssymb,amscd,amsthm,euscript}
\usepackage[utf8]{inputenc}
\usepackage[T2A]{fontenc}
\usepackage[russian]{babel}
\usepackage[metapost]{mfpic}
\usepackage{epstopdf}
\usepackage{fancybox,textcomp}
\usepackage{calrsfs}
\usepackage{enumerate,comment}
\usepackage[margin=2cm]{geometry}
\usepackage{secdot}

\theoremstyle{plain}
\newtheorem{thm}{\bf Теорема}[section]

\theoremstyle{remark}
\newtheorem{dfn}{\bf Определение}[section]

\newcommand*{\set}[1]{\left\{#1\right\}}

\def\lc{\ulcorner}
\def\rc{\urcorner}
\def\a{\alpha}
\def\dl{\delta}
\def\Dl{\Delta}
\def\lm{\lambda}

\def\vp{\varphi}
\def\<{\langle}
\def\>{\rangle}

\def\ul{\underline}

\def\eq{\leftrightarrow}
\def\eqq{\longleftrightarrow}

\def\sbs{\subset}

\def\0{\varnothing}
\def\опр{\overset\text{Опр}\to\eqq}
\def\les{\leqslant}
\def\ges{\geqslant}

\def\R{\mathbb R}
\def\Pcal{\mathcal P}
\def\Scal{\mathcal S}
\def\Tcal{\mathcal T}

\allowdisplaybreaks

\begin{document}


\noindent
{\bf О понятии рода структуры (в смысле Бурбаки)}
\\
\\
\textit{А. Х. Назиев}
\bigskip

\section*{\bf Введение}
При изучении различных математических объектов~--- групп, колец,
полей, векторных пространств и т. д.~--- обычно подчёркивают,
что природа элементов, из которых состоят эти объекты, совершенно
несущественна, а также отмечают, что изоморфные объекты в математике
принято не различать, считая их как бы различными копиями одного
и того же объекта. Иногда даже эти два положения называют
отличительными чертами всей современной математики.

В вузовских курсах алгебры, геометрии, математического анализа
указанные положения сообщаются, как правило, в виде своего рода
психологических установок, а не в виде правил, имеющих точный
математический смысл. И это понятно, ибо точный смысл названные
положения приобретают не в рамках конкретных теорий, изучающих
те или иные математические объекты, а в рамках некоторой специальной
теории~--- общей теории родов структуры. Началам этой теории и посвящена
настоящая работа.

\smallskip
Могло бы показаться, что для того, чтобы придать точный
математический смысл положению о несущественности природы элементов,
из которых состоят математические объекты, нужно сначала дать точное
определение понятию природы элементов. Однако представляется по меньшей
мере неэкономным трудиться над определением чего-нибудь только для того,
чтобы затем сказать: <<Вот именно это для нас и несущественно>>.
Поэтому не удивительно, что поступают иначе.

Если совершенно несущественно, из каких именно элементов устроен данный
объект, то любой его элемент может быть заменён любым другим предметом,
лишь бы этот предмет был отличен от всех других элементов, составляющих
рассматриваемый объект. Ясно, что посредством таких замен можно перейти
от данного множества к любому другому множеству, содержащему <<столько
же>>\ элементов, сколько данное. <<Столько же>>\ на математическом
языке означает существование биекции данного множества на это
<<другое>>\ множество. Таким образом, несущественность природы
элементов, составляющих математические объекты, означает, что, если
на данном множестве имеется некоторая структура, и существует биекция
этого множества на другое множество, то и на этом другом множестве
имеется <<такая же>>\ структура.

Описанное только что свойство называется переносимостью и является
основным в теории структур (построенной Н. Бурбаки). При этом оно
формулируется таким образом, что из определения тотчас же вытекает,
что изоморфные объекты имеют одинаковые теории, так что тут же получает
точный смысл и второе из упомянутых в начале положений.

\medskip
В последние десятилетия не раз предпринимались попытки донести теорию структур
Бурбаки до широкой аудитории. Все эти попытки представляются нам неудачными.
И дело здесь не в том, что было много недостатков в сказанном. Дело в том,
что при этом ничего не говорилось о чрезвычайно важном понятии переносимости.
Сказать, что такое структура, совсем не трудно\footnote{Хотя и здесь не обошлось
без курьёзов: в нескольких работах под названием структуры (и со ссылкой на
Бурбаки!) определяется не понятие структуры, а понятие алгебраической системы.
Одна таких работ защищена недавно в качестве докторской диссертации по методике
преподавания математики. Особенно забавно то, что в этой работе её автор пишет:
<<``структуры являются орудиями математики'', и только через них можно в определённой
степени систематизировать математику, дать общее представление о ней>>. Получается,
что только через алгебраические системы можно в определённой степени систематизировать
математику, дать общее представление о ней! А что же топология, теория вероятностей,
математическая статистика и другие не алгебраические разделы математики?}.
Однако этого не достаточно. Нужно и важно ещё (и сделать это гораздо труднее)
объяснить, что означают слова <<такая же структура>>. Это и делается с помощью
понятия переносимости, что приводит в итоге к понятию {\em рода} структуры.

С другой стороны, теория структур Бурбаки в России (СССР) не раз подвергалась нападкам.
Понимание структуры Н.~Бурбаки называли ограниченным и односторонним. Предлагали его
расширить. Например, высказывалось (весьма авторитетным математиком) пожелание включить
в понятие структуры математические модели реальных явлений. Это пожелание представляется
нам по меньшей мере неясным. В самом деле, разве не являются математическими моделями
реальных явлений элементарные функции, дифференциальные и интегральные уравнения,
топологические векторные (и, в частности, гильбертовы) пространства, линейные операторы
и их спектры, банаховы алгебры, меры (в частности, мера Винера)? А ведь всё это излагается
в трактате Бурбаки.

Другой упрёк в адрес Бурбаки, высказанный одним известным методистом, состоит в том, что
теорией Бурбаки, якобы, не охватываются комбинаторные структуры. Это просто неверно.
Элементарная комбинаторика превосходно изложена самим Бурбаки в первом томе его трактата,
высшие разделы~--- в известной монографии Э.~Баннаи и Т.~Ито, написанной на том же языке,
что и ``Спектральные теории'' Бурбаки.

По поводу  этих и подобных им попыток необходимо заметить следующее. Желающие указать
примеры структур, не охватываемых теорией Бурбаки, должны представлять себе всю грандиозность
проблемы, которую они берутся решать. Бурбаки выделил не ``три структуры'', как думают
и говорят многие, а три базисных {\em разновидности родов\/} структуры. Сам же Бурбаки
и допустил, что, может быть, когда-нибудь будут обнаружены роды структуры, не сводимые
к родам этих трёх разновидностей. Но пока этого не произошло. И вряд ли скоро произойдёт.
Чтобы понять --- почему, достаточно вспомнить историю попыток решения в радикалах
уравнений пятой степени. Легко ли было установить, что не все уравнения пятой степени
разрешимы в радикалах? Даже сейчас нет простого доказательства этого факта. А в случае
со структурами нужно доказать~--- именно доказать, ибо теория структур Бурбаки является
математической теорией~--- существование рода структуры, не сводимого к родам трёх
описанных Бурбаки разновидностей. Задача --- неизмеримо более сложная, если вообще
разрешимая.

\medskip
Ещё одно (завуалированное) направление атак на теорию Бурбаки состоит в том,
что предлагают вместо неё что-нибудь, по мнению предлагающих, более подходящее.
Например, предлагают считать математику наукой о математических моделях, или о
видах математических моделей, или о математических схемах. Эти авторы, при всей
их высокой авторитетности, не замечают (или не хотят замечать) одного весьма
существенного обстоятельства: они высказывают свои предложения на \emph{эвристическом}
уровне, тогда как Бурбаки построил \emph{\bfseries математическую теорию} математических
структур. И нет никакого сомнения: если бы эти авторы взяли на себя труд реализовать
свои предложения в виде \emph{\bfseries математических теорий} (например, теории видов
математических моделей), то получилось бы у них практически то же самое,
что у Бурбаки, только с другими названиями.
$$
  *\quad*\quad*
$$

Изложение теории родов структуры у Н.~Бурбаки хотя и занимает всего несколько
страниц~\cite[С. 246--254]{bour}, всё же является достаточно трудным.
Однако в упражнениях Н.~Бурбаки намечает другой подход, основанный не на
понятии схемы конструкции ступени, как в основном тексте, а на понятии
типа ступени. При этом, полагаясь в допустимых пределах на здравый смысл,
можно обойтись без усложняющего рассмотрения понятия равновесного
знакосочетания, и тогда получится вполне доступное и в то же время
достаточно строгое изложение. Именно так мы и поступаем, следуя~\cite{AN:dissB}.
При этом, чтобы сделать изложение как можно более доступным, мы включили
в него большое количество примеров.

\smallskip
Работа заканчивается определением понятия рода структуры. Автор надеется продолжить
работу, включив в неё рассмотрение трёх основных видов родов структуры, комбинированных
структур и структур с морфизмами.
$$
  *\quad*\quad*
$$

\paragraph*{В работе используются следующие обозначения:}
$\Scal$ --- класс всех множеств;
$\Dl$ --- тождественная функция на~$\Scal$;
$\Pcal(X)$ --- множество всех подмножеств множества $X$;\
$I_X$ --- тождественная функция на множестве $X$. Если ясно или неважно,
о каком именно $X$ идёт речь, пишем просто $I$.

\smallskip
Остальные обозначения вводятся по ходу дела.

\section{\bf Ступени и структуры}

\begin{dfn}
Пусть $X_1,\ \dots,\ X_n$~--- множества. Понятие {\it ступени над множествами}
$X_1,\ \dots,\ X_n$ вводится следующими соглашениями.
\begin{enumerate}[$(S1)$]
\item Каждое из множеств $X_1,\ \dots,\ X_n$~--- ступень над множествами
      $X_1,\ \dots,\ X_n$.
\item Если $X$~--- ступень над множествами $X_1,\ \dots,\ X_n$, то и
      $\Pcal(X)$~--- ступень над $X_1,\ \dots,\ X_n$.
\item Если $X$, $Y$~--- ступени над $X_1,\ \dots,\ X_n$, то и $X\times Y$
     ~--- ступень над $X_1,\ \dots,\ X_n$.
\item Других ступеней над $X_1,\ \dots,\ X_n$ нет.
\end{enumerate}

Совокупность всех ступеней над $X_1$, \dots, $X_n$ называется {\it ле\-стницей
множеств}, построенной над $X_1$, \dots, $X_n$, а элементы ступеней~---
{\it структурами} над $X_1$, \dots, $X_n$.
\end{dfn}

\subsubsection*{Примеры.} 1. Любое из множеств
$$
\begin{array}{l}
 X_i,\ X_i\times X_j,\ X_i\times(X_j\times X_k),\\
 \Pcal(X_i),\ \Pcal(X_i\times X_j),\ \Pcal(X_i\times(X_j\times X_k)),\\
 X_j\times\Pcal(X_j\times X_k),\ \Pcal(X_i\times X_j)\times \Pcal(X_i\times(X_j\times X_k))
\end{array}
$$
есть ступень над $X_1,\ \dots,\ X_n$.
Но любая ступень $S$ над $X_1,\ \dots,\ X_n$ является элементом ступени
$\Pcal(S)$ и потому~--- структурой над $X_1,\ \dots,$ $X_n$. Значит, все
указанные множества являются и структурами над $X_1,\ \dots,\ X_n$.

\smallskip
2. Как отмечено в примере 1, любое из множеств $X_1$, \dots, $X_n$~---
ступень над $X_1$, \ldots, $X_n$. Значит, любой элемент любого из множеств
$X_1$, \ldots, $X_n$~--- структура на этих множествах. Такие структуры
называют {\it выделенными элементами}.

\smallskip
3. ($n = 1$) Отношение порядка на $\R$ есть структура на этом множестве.
Действительно, это отношение является подмножеством в $\R\times\R$
и потому~--- элементом ступени $\Pcal(\R\times\R)$. Вообще, любое
отношение на любом множестве есть структура на этом множестве.

\smallskip
4. ($n = 1$) Операция сложения на $\R$ есть структура на этом
множестве. Действительно, эта операция есть функция на $\R\times\R$
со значениями в $\R$, значит~--- подмножество в $(\R\times\R)\times\R$
и потому~--- элемент ступени $\Pcal((\R\times\R)\times\R)$. Вообще,
любая операция на любом множестве есть структура на этом множестве.

\smallskip
5. ($n = 1$) Топология на $\R$ есть структура на этом множестве.
Действительно, топология на $\R$ есть совокупность всех открытых
подмножеств в $\R$, значит,~--- подмножество множества $\Pcal(\R)$
и потому~--- элемент ступени $\Pcal(\Pcal(\R))$. Вообще, любая
топология на любом множестве есть структура на этом множестве.

\smallskip
6. ($n = 1$) Пусть $\omega$~--- бинарная операция, $\rho$~---
бинарное отношение на множестве $X$. Тогда $\omega\in\Pcal((X\times X)
\times X)$, $\rho\in\Pcal(\Pcal(X))$, значит $(\omega,\ \rho)\in
\Pcal((X\times X)\times X)\times\Pcal(\Pcal(X))$. Тем самым, $(\omega,\ \rho)$~---
структура на $X$. Вообще, любая упорядоченная пара структур
на любом множестве есть структура на этом множестве. Чаще всего
встречаются упорядоченные пары, составленные из бинарной операции
и отношения порядка, бинарной операции и топологии, бинарного отношения
и топологии.

\smallskip
7. ($n = 1$) Пусть $\omega$~--- бинарная операция, $\rho$~---
бинарное отношение, $\tau$~--- топология на множестве $X$.
Тогда упорядоченная тройка $(\omega,\ \rho,\ \tau)$~--- структура на $X$.
Вообще, любой упорядоченный набор структур на любых множествах
$X_1,\ \dots,\ X_n$ есть структура на этих множествах.

\smallskip
8. ($n = 2$) Пусть $V$~--- векторное пространство над полем $k$, $+$~---
операция сложения векторов, $\cdot$~--- операция умножения их на скаляры.
Тогда $(+,\ \cdot)$~--- структура на множествах $V,\ k$. По некоторым
причинам, которые выяснятся позже, говорят также иначе: $(+,\ \cdot)$
есть структура на $V$ при вспомогательном множестве $k$.

\smallskip
9. ($n = 2$) Пусть $d$~--- метрика на множестве $X$, то есть
функция на $X\times X$ со значениями в $\R$, при всех $x$, $y$,
$z\in X$ удовлетворяющая условиям:  $d(x,\; y) \ges 0$;
$d(x,\; y) = 0 \eq x = y$; $d(x,\; y) = d(y,\; x)$;
$d(x,\; z) \les d(x,\; y) + d(y,\; z)$.
Тогда $d$~--- структура на $X$ и $\R$. И в этом случае говорят также, что $d$
есть структура на $X$ при вспомогательном множестве~$\R$.

\medskip
Пусть $X_1,\ X_2$\/ и $Y_1,\ Y_2$~--- множества. Рассмотрим
ступени: $\Pcal(X_1\times X_2)$~--- над $X_1,\ X_2$\/ и
$\Pcal(Y_1\times Y_2)$~--- над $Y_1,\ Y_2$. Это, вообще говоря,~---
различные ступени над различными множествами. Тем не менее мы
хорошо осознаём, что эти ступени построены <<по одному
способу>>, или, как ещё говорят в таких случаях, эти ступени
имеют один и тот же тип. Наша ближайшая цель~--- дать точное
определение понятию типа ступени. Но сначала~--- несколько
предварительных определений.

Один из способов построения ступени над данными $n$ множествами
состоит в выборе одного из них: по данным $X_1,\ \dots,$ $X_n$ выбираем
$X_i$, где $i$~--- одно из чисел $1,\ \dots,\ n$. Обозначим эту операцию
через $pr^n_i$, то есть условимся, что $pr^n_i$ есть функция на $\Scal^n$
со значениями в $\Scal$, действующая по правилу: $pr^n_i(X_1,\ \dots,\ X_n)
= X_i$, $i = 1,\ \dots,\ n$. Если ясно, о каком $n$ идёт речь, то
верхний индекс не пишут.

Другой способ построения ступени состоит в том, что, располагая уже готовой
ступенью $X$, образуют ступень $\Pcal(X)$. Операция образования множества
$\Pcal(X)$ по данному множеству $X$ есть функция $\Pcal: \Scal\to\Scal$.

Наконец, ещё один способ построения ступени состоит в том, что,
располагая уже готовыми ступенями $S_1,\ S_2$, образуют ступень
$S_1\times S_2.$ Операция перехода от данных $X_1\ \dots,\ X_n$
к $S_1\times S_2$ оказывается произведением операций перехода
от данных $X_1\ \dots,\ X_n$ к $S_1$ и от них же~--- к $S_2$;
произведением~--- в смысле следующего определения.

\begin{dfn}
Пусть $F_1$ и $F_2$~--- функции на $\Scal^n$ со значениями в $\Scal$.
{\em Произведением $F_1$ и $F_2$} называется функция
$$
  F_1\otimes F_2: \Scal^n \to \Scal,
$$
действующая по правилу:
$$
  (F_1\otimes F_2)(X_1,\ \dots,\ X_n) = F_1(X_1,\ \dots,\ X_n)
                                    \times F_2(X_1,\ \dots,\ X_n).
$$
\end{dfn}

Теперь мы готовы дать определение понятию типа ступени. В силу этого
определения каждый тип ступени над $n$ множествами оказывается
некоторым отображением из $\Scal^n$ в $\Scal,$ но, конечно, не все такие
отображения являются типами ступени. О том, какие именно являются,
и говорится в следующем определении.

\begin{dfn}
Понятие {\it типа ступени над $n$ множествами} вводится следующими
соглашениями.
\begin{enumerate}
\item Каждое из отображений $p_j^n: \Scal^n\to \Scal$, $i = 1,\ \dots,\ n,$
      является типом ступени над $n$ множествами.
\item Если $T$~--- тип ступени над $n$ множествами, то и $\Pcal\circ T$
     ~--- тоже.
\item Если $S$, $T$~--- типы ступени над $n$ множествами, то и $S\otimes T$
     ~--- тоже.
\item Других типов ступени над $n$ множествами нет.
\end{enumerate}
\end{dfn}

\begin{dfn}
Пусть $T$~--- тип ступени над $n$ множествами и $X_1$, \dots, $X_n\in\Scal$.
Множество $T(X_1, \ldots, X_n)$~--- значение функции $T$ в точке
\mbox{$(X_1, \ldots, X_n)$}~--- называют {\it реализацией типа ступени
$T$ на множествах} $X_1$, \ldots, $X_n$. Ясно, что это~--- ступень над
указанными множествами. Её элементы называют {\it структурамиютипа} $T$
на множествах $X_1$, \ldots, $X_n$.
\end{dfn}

\subsubsection*{Примеры.} I. {\it Типы ступени над одним множеством.}

1. \ul{$T = pr_1^1.$} Это~--- $\Delta.$ Действует по правилу: $T(X) = X.$
Структуры этого типа~--- выделенные элементы (данного множества).

\smallskip
2. \ul{$T = \Pcal$.}\ \ Является типом ступени (над одним множеством), потому
что им является $\Delta,$ а $\Pcal = \Pcal\circ\Delta.$ Действует по правилу:
$T(X) = \Pcal(X).$ Структуры этого типа~--- подмножества основного
множества.

\smallskip
3. \ul{$T = \Delta\otimes\Delta.$} Действует по правилу: $T(X) = X\times X.$
Структуры этого типа~--- упорядоченные пары элементов основного множества.

\smallskip
4. \ul{$T = \Delta\otimes \Pcal.$} Действует по правилу:
$$
  T(X) = X\times \Pcal(X).
$$
Структуры этого типа~--- упорядоченные пары, образованные элементом
и подмножеством основного множества.

\smallskip
5. \ul{$T = \Pcal\circ\Pcal.$} Действует по правилу:
$$
  T(X) = \Pcal(\Pcal(X)).
$$
Структуры этого типа~--- совокупности подмножеств основного множества.
Таковы, например, топологии.

\smallskip
6.  \ul{$T = \Pcal\circ(\Delta\otimes\Delta).$} Действует по правилу:
$$
  T(X) = \Pcal(X\times X).
$$ Структуры этого типа~--- подмножества второй
декартовой степени основного множества. Таковы, например, бинарные
отношения на основном множестве.

\smallskip
7. \ul{$T = \Pcal\circ((\Delta\otimes\Delta)\otimes\Delta).$} Действует
по правилу:
$$
  T(X) = \Pcal((X\times X)\times X).
$$
Структуры этого типа~--- тернарные отношения на основном множестве. Таковы, например,
бинарные операции.

\smallskip
8. \ul{$T = (\Pcal\circ \Pcal)\otimes(\Pcal\circ((\Delta\otimes\Delta)\otimes\Delta)).$}
Действует по правилу:
$$
  T(X) = \Pcal(\Pcal(X))\times\Pcal((X\times X)\times X).
$$
Структурами этого типа являются, например, упорядоченные пары, образованные
топологиями и бинарными операциями.

\medskip
II. {\it Типы ступени над двумя множествами}.

\smallskip
1.  \ul{$T = pr^2_j.$} Действует по правилу: $T(X_1,\ X_2) = X_j.$
Структуры этого типа~--- выделенные элементы: при $j = 1$~--- первого
множества, при $j = 2$~--- второго множества.

\smallskip
2. \ul{$T = pr^2_i\otimes pr^2_j.$} Действует по правилу:
$T(X_1,\ X_2) = X_i\times X_j.$ Например, при $i = 1,$ $j = 2$,~---
по правилу $T(X_1,\ X_2) = X_1\times X_2.$ Структуры этого типа~---
выделенные элементы соответствующего декартова произведения.

\smallskip
3. \ul{$T = \Pcal\circ pr^2_j.$} Действует по правилу:
$T(X_1,\ X_2) = \Pcal(X_j),$ $j = 1,\ 2.$ Структуры этого типа~---
подмножества одного из данных множеств: первого, если $j = 1;$
второго, если $j = 2.$

\smallskip
4. \ul{$T = pr^2_i\otimes(\Pcal\circ pr^2_j).$} Действует по правилу:
$$
  T(X_1,\ X_2) = X_i\times\Pcal(X_j),
$$
$j = 1,\ 2.$ Структуры этого типа~--- упорядоченные пары, первые
члены которых~--- элементы $X_i$, вторые~--- подмножества $X_j.$

\smallskip
5. \ul{$T = \Pcal\circ(pr^2_i\otimes pr^2_j).$} Действует по правилу:
$$
  T(X_1,\ X_2) = \Pcal(X_i\times X_j),
$$
$j = 1,\ 2.$ Структуры этого типа~--- подмножества произведения
$X_i\times X_j$, то есть соответствия между $X_i$ и $X_j.$

\smallskip
6. \ul{$T = \Pcal\circ((pr^2_i\otimes pr^2_j)\otimes pr^2_k)$.}
Пусть, для определён\-ности, $i = j = 1,$ $k = 2.$
Тогда $T$ действует по правилу:
$$
  T(X_1,\ X_2) = \Pcal((X_1\times X_1)\times X_2).
$$

Структуры этого типа~--- подмножества декартова произведения
квадрата первого множества на второе. Таковы, в частности, отображения
из квадрата первого множества во второе. Например, метрика
на множестве $S$ есть отображение из $S\times S$ в $\R$.

\section{\bf Канонические распространения функций}

\begin{dfn}
Пусть $f: X\to Y.$ Определим $\Pcal\langle f\rangle: \Pcal(X)\to \Pcal(Y)$
условием: $\Pcal\langle f\rangle(s) = f[s].$ Так определённая функция
называется {\it каноническим распространением функции $f$ на подмножества}.
\end{dfn}

\begin{thm}[Свойства канонических распространений функций на подмножества]\hfill
\begin{enumerate}[$1^\circ$]
\item $\Pcal\langle \Delta_X\rangle = \Delta_{\Pcal(X)}.$

\item $\Pcal\langle g\circ f\rangle = \Pcal\langle g \rangle \circ
                                      \Pcal\langle f \rangle.$

\item Если $f$ инъективна (сюръективна, биективна), то и
          $\Pcal\langle f \rangle$~--- тоже.

\item Если $f$ обратима, то и $\Pcal\langle f \rangle$~--- тоже,
          причём
          $$
            (\Pcal\langle f \rangle)^{-1} = \Pcal\langle f^{-1} \rangle.
          $$
\end{enumerate}
\end{thm}

\begin{dfn}
Теперь пусть $f_1: X_1\to Y_1,$ $f_2: X_2\to Y_2.$ Определим функцию
$$
  f_1 * f_2: X_1\times X_2 \to Y_1\times Y_2
$$
условием:
$$
  (f_1 * f_2)(x_1,\ x_2) = (f_1(x_1),\ f_2(x_2)).
$$
Так определённая функция называется {\it каноническим распространением
функций $f_1$ и~$f_2$ на произведение}.
\end{dfn}

\begin{thm}[Свойства канонических распространений функций
             на произведения]\hfill
\begin{enumerate}[$1^\circ$]
\item $\Delta_{X_1} * \Delta_{X_2} = \Delta_{X_1\times X_2}.$
\item $(g_1\circ f_1)*(g_2\circ f_2) = (g_1 * g_2)\circ(f_1 * f_2).$
\item Если $f_1,$ $f_2$ инъективны (сюръективны, биективны), то и
          $f_1 * f_2$~--- тоже.
\item Если $f_1,$ $f_2$ обратимы, то и $f_1 * f_2$~--- тоже,
          причём
          $$
            (f_1 * f_2)^{-1} = f_1^{-1} * f_2^{-1}.
          $$
\end{enumerate}
\end{thm}

\begin{dfn}
Наконец, пусть $T$~--- произвольный тип ступени над $n$ множествами
($n \ges 1$) и $f_i: X_i\to Y_i,$ $i = 1,\ \dots,\ n.$ Определим функцию
$$
  T\langle f_1,\ \dots,\ f_n\rangle : T(X_1,\ \dots,\ X_n)\to T(Y_1,\ \dots,\ Y_n)
$$
условием:
$$
  T\langle f_1,\ \dots,\ f_n\rangle =
   \left\{
    \begin{array}{rll}
     f_i, & \mbox{если} & T = pr^n_i;\\
     \Pcal\langle T_1\langle f_1,\ \dots,\ f_n\rangle\rangle,
         & \mbox{если} & T = \Pcal\circ T_1;\\
     T_1\langle f_1,\ \dots,\ f_n\rangle * T_2\langle f_1,\ \dots,\ f_n\rangle,
         & \mbox{если} & T = T_1 \otimes T_2.
    \end{array}
   \right.
$$
Так определённая функция $T\langle f_1,\ \dots,\ f_n\rangle$ называется
{\em каноническим распространением функций $f_1,\ \dots,\ f_n$ по типу~$T$}.
\end{dfn}

\begin{thm}[Основные свойства канонических распространений функций]\hfill
\begin{enumerate}[$1^\circ$]
\item $T\langle I_{X_1},\ \dots,\ I_{X_n}\rangle =
             I_{T\langle X_1\ \dots,\ X_n\rangle}$;
\item $T\langle g_1\circ f_1,\ \dots,\ g_n\circ f_n\rangle =
             T\langle g_1,\ \dots,\ g_n\rangle \circ
             T\langle f_1,\ \dots,\ f_n\rangle;$
\item Если $f_1,\ \dots,\ f_n$ инъективны (сюръективны,
            биективны), то и \newline $T\langle f_1\ \dots,\ f_n\rangle$~--- тоже;
\item Если $f_1,\ \dots,\ f_n$ обратимы,
            то и $T\langle f_1\ \dots,\ f_n\rangle$~--- тоже,
            причём
            $$
              (T\langle f_1\ \dots,\ f_n\rangle)^{-1} =
               T\langle f_1^{-1}\ \dots,\ f_n^{-1}\rangle.\qquad\qquad\qquad
            $$
\end{enumerate}
\end{thm}

Несложные доказательства теорем Т3.1--Т3.3 оставляются читателю
в качестве упражнений.

\subsubsection*{Замечание.} Хотелось бы писать `$\Pcal(f)$' вместо
`$\Pcal\langle f \rangle$', `$S\times T$' вместо `$S\otimes T$'
и `$T(f_1,\ \dots,\ f_n)$' вместо `$T\langle f_1,\ \dots,\ f_n\rangle$'.
Мы вынуждены отказаться от этого, потому что и `$\Pcal(f)$', и `$S\times T$',
и `$T(f_1,\ \dots,\ f_n)$' уже были определены ранее: первое~--- как степень~$f$,
второе~--- как декартово произведение $S$ и $T$, третье~--- как реализация
типа ступени $T$ на множествах $f_1,\ \dots,\ f_n.$

\subsubsection*{Примеры.} 1. \ul{$n = 1,\ T = \Delta.$} Пусть $X$, $X' \in \Scal$,
$f: X\to X'.$ Тогда
$$\begin{array}{l}
   T(X) = X,\ T(X') = X',\\
   T\< f\> = \Dl\< f \> = pr^1_1\<f\> = f.
  \end{array}
$$
Таким образом, $T\<f\>$ каждой структуре $s\in T(X) = X$
ставит в соответствие структуру $s'\in T\<f\>(s) = f(s)\in T(X') = X'.$
Например, если $X = Y = \R,$ $f = \lm_x(x-3)$ и $s = 1$, то $s' = -2.$

\smallskip
2. \ul{$n = 1,\ T = \Pcal.$} Пусть $X$, $X' \in \Scal$, $f: X\to X'.$
Тогда
$$\begin{array}{l}
   T(X) = \Pcal(X),\ T(X') = \Pcal(X'),\\
   T\< f\> = \Pcal\<f\>: \Pcal(X)\to \Pcal(X').
  \end{array}
$$
Для любого $s\in T(X) = \Pcal(X)$ имеем $T\<f\>(s) = p\<f\>(s)
= f[s].$ Таким образом, в данном случае $T\<f\>$ каждой структуре
$s\in T(X)$ ставит в соответствие образ $s$ относительно~$f$.
Пусть, например, $X = X' = \R$, $f = \lm_x|x|$ и $s = [-3; 2].$
Тогда
$$
  T\<f\>(s) = f[s] = \{|x|: x \in [-3; 2]\} = [0; 3].
$$

\smallskip
3. \ul{$n = 1,\ T = \Dl\otimes\Dl$.} Пусть $X$, $X' \in \Scal$, $f: X\to X'.$
Тогда
$$\begin{array}{l}
   T(X) = X\times X,\ T(X') = X'\times X',\\
   T\< f\> = (\Dl\otimes\Dl)\<f\> = \Dl\<f\> * \Dl\<f\> = f * f.
  \end{array}
$$
Каждой $s = (s_1,\ s_2) \in T(X) =
X\times X$ отображение $T\<f\>$ ставит в соответствие
$$
  s' = T\<f\>(s) = (f*f)(s) = (f*f)(s_1,\ s_2) = (f(s_1),\ f(s_2)).
$$
Например, если $X=X'=\R$, $f = \lm_x|x|$ и $s = (1; -4)$, то $s' = (1,\ 4).$

\smallskip
4. \ul{$n = 1,\ T = \Pcal\circ \Pcal$.} Для любых $X,\ X' \in \Scal$ и $f: X\to X'$
имеем:
$$\begin{array}{l}
   T(X) = \Pcal(\Pcal(X)),\ T(X') = \Pcal(\Pcal(X')),\\
   T\<f\>: \Pcal(\Pcal(X))\to \Pcal(\Pcal(X')).
  \end{array}
$$
Каждой $s\in T(X)$ отображение $T\<f\>$ ставит в соответствие
$$\begin{array}{rcl}
   s' &=& T\<f\>(s) = \Pcal\<\Pcal\<\<f\>\> = \Pcal\<f\>[s]\\
      &=& \{\Pcal\<f\>(u): u\in s\} = \{f[u]: u\in s\}.
  \end{array}
$$
Таким образом, $T\<f\>(s)$ в этом случае состоит из всех образов
множеств, принадлежащих $s$, относительно $f$. Если $X=X'=\R$,
$s$~--- обычная топология на $\R$ и $f = \lm_x(x-3)$, то $s'$~---
это снова та же топология, ибо элементы структуры $s'$ получаются
из элементов структуры $s$ сдвигом на 3 единицы влево, а топология
прямой, как известно, инвариантна относительно сдвигов. А вот
для $f = \lm_x|x|$ получим $s'$, отличную от $s$. Более того,
элементы $s'$ уже не обязательно будут открытыми подмножествами
прямой. Например, если $u = (-1; 1)$, то $u$ открыто, но $f[u]
= [0; 1)$~--- не открыто.

\smallskip
5. \ul{$n = 1,\ T = \Pcal\circ(\Dl\otimes\Dl)$.} Для любых $X$, $X' \in \Scal$,
$f: X\to X'$ имеем:
$$\begin{array}{l}
   T(X) = \Pcal(X\times X),\ T(X') = \Pcal(X'\times X'),\\
   T\< f\> = \Pcal\<(\Dl\otimes\Dl)\<f\>\>
           = \Pcal\<\Dl\<f\> * \Dl\<f\>\> = \Pcal\<f * f\>.
  \end{array}
$$
Каждой $s \in T(X) = \Pcal(X\times X)$ отображение $T\<f\>$ ставит в соответствие
$$\begin{array}{rcl}
   s' &=& T\<f\>(s) = \Pcal\<(f*f)\>(s) = (f*f)[s]\\
      &=& \{(f(x),\ f(y)): (x,\ y) \in s\}.
  \end{array}
$$

Отсюда видно, что
  \begin{align*}
  (x',\ y') \in s' &\leftrightarrow (\exists x)(\exists y)(x' = f(x) \land
                                                y' = f(y) \land
                                                (x,\ y) \in s)    \\
                   &\leftrightarrow (\exists x)(\exists y)((x,\ x') \in f \land
                                                (y,\ y') \in f \land
                                                (x,\ y) \in s)    \\
                   &\leftrightarrow (\exists x)(\exists y)((x',\ x) \in f^{-1} \land
                                                (x,\ y) \in s  \land
                                                (y,\ y') \in f)   \\
                   &\leftrightarrow (\exists x)((x',\ x) \in f^{-1} \land
                                     (\exists y)((x,\ y) \in s  \land
                                                 (y,\ y') \in f)) \\
                   &\leftrightarrow (\exists x)((x',\ x) \in f^{-1} \land
                                     (x,\ y') \in f\circ s)       \\
                   &\leftrightarrow (x',\ y') \in f\circ s\circ f^{-1}.
   \end{align*}
Таким образом, $s' = f\circ s\circ f^{-1}.$ Эта связь между $s$ и $s'$
становится особенно прозрачной, если $s$ и $f^{-1}$~--- функции (т. е.
$s$~--- функция и $f$ --- обратимая функция). Тогда эта связь может
быть представлена в виде следующей коммутативной диаграммы:
$$
  \begin{CD}
       X       @>s>>\!\!\! X        \\
   @Af^{-1}AA       \!\!\! @VVfV      \\
       X'      @>s'>> \  X'\ .
  \end{CD}
$$
(Коммутативность этой диаграммы как раз и состоит в том, что
$s' = f\circ s\circ f^{-1}.$)

Кроме того, если $f$~--- биекция, то, независимо от того, является
$s$ функцией или нет:
$$\begin{array}{rcl}
   (x',\ y') \in s' &\eq& (\exists x,\ y)(x' = f(x) \land
                                        y' = f(y) \land
                                       (x,\ y) \in s)    \\
                   &\eq& (\exists x,\ y)(x = f^{-1}(x') \land
                                        y = f^{-1}(y') \land
                                       (x,\ y) \in s)    \\
                   &\eq& (f^{-1}(x'),\ f^{-1}(y')) \in s.
   \end{array}
$$

Изображая $T\< f\>(s)$ на чертеже, полезно помнить п. $2^\circ$
T3.2. В силу этого предложения $f*f = (I*f)\circ(f*I).$ Поэтому
переход от $s$ к $s'$ можно осуществлять в два этапа: сначала найти
образ $s$ относительно $f*I,$ а затем~--- образ полученного
множества относительно $I*f.$ Для $X = X' = \R$ и $f = \lm_x|x|$,
например, это будет означать, что сначала отражаются относительно
горизонтальной оси те части $s$, которые лежат ниже этой оси,
а затем отражаются относительно вертикальной оси те части
полученного множества, которые лежат левее этой оси. Если
взять в качестве $s$ отношение нестрого порядка на $\R$, то
в результате этих действий получится первый координатный угол
(и, значит,~--- не отношение порядка). А если взять в качестве
$s$ параболу, определяемую уравнением $y = x^2 + 2x - 3,$
график которой изображён на рис. а) ниже, то в результате
этих действий получится фигура, изображённая на рис. б).

\begin{center}
\includegraphics{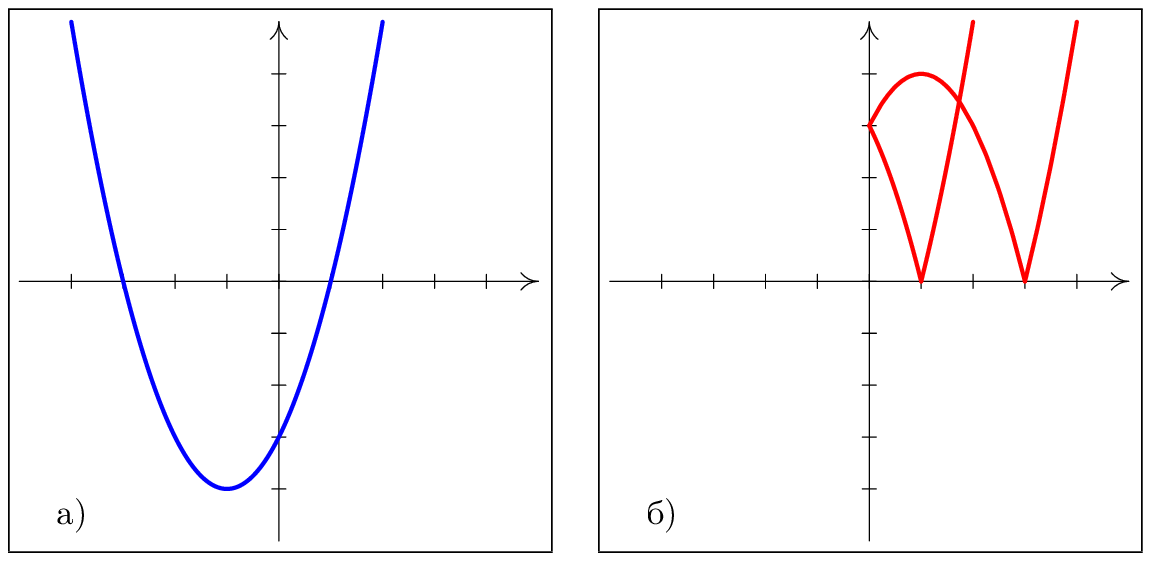}
\end{center}

\medskip
6. \ul{$n = 1,\ T = \Pcal\circ((\Dl\otimes\Dl)\otimes\Dl)$.}
Для любых $X$, $X' \in \Scal$, $f: X\to X'$ имеем:
$$\begin{array}{l}
   T(X) = \Pcal((X\times X)\times X),\ T(X') = \Pcal((X'\times X')\times X'),\\
   T\< f\> = \Pcal\<((\Dl\otimes\Dl)\otimes\Dl)\<f\>\>
           = \Pcal\<(f * f) * f\>.
  \end{array}
$$
Каждой $s \in T(X)$ отображение $T\<f\>$ ставит в соответствие
$$
\begin{array}{rcl}
   s' &=& T\<f\>(s) = \Pcal\<((f*f)*f)\>(s) = ((f*f)*f)[s]\\
      &=& \{((f(x),\ f(y)),\ f(z)): ((x,\ y),\ z) \in s\}.
  \end{array}
$$

Отсюда видно, что
\begin{align*}
  ((x',\ y'),\ z') \in s' &\eq (\exists x, y, z)(x' = f(x) \land
                                                  y' = f(y) \land
                                                  z' = f(z) \land
                                                ((x,\ y),\ z) \in s)    \\
&\eq  (\exists x, y, z)((x,\ x')\in f   \land
                                                  (y,\ y')\in f   \land
                                                  (z,\ z')\in f   \land
                                                ((x,\ y),\ z) \in s)    \\
                         &\eq (\exists x, y, z)((x',\ x)\in f^{-1}   \land
                                                  (y',\ y)\in f^{-1}   \land
                                                  ((x,\ y),\ z) \in s   \land
                                                  (z,\ z')\in f)       \\
                         &\eq (\exists x, y, z)(((x',\ y'),\ (x,\ y))
                                                   \in f^{-1}*f^{-1}   \land
                                                    ((x,\ y),\ z) \in s  \land
                                                   (z,\ z')\in f)      \\
                         &\eq (\exists x, y, z)(((x',\ y'),\ z')\in
                                                  f\circ s\circ
                                                  (f^{-1}*f^{-1})).
\end{align*}

Таким образом, $s' = f\circ s\circ(f^{-1}*f^{-1}).$ Если $s$ и $f^{-1}$~---
функции (т.~е. $s$~--- функция и $f$~--- обратимая функция), то эта
связь между $s$ и $s'$ представляется коммутативной диаграммой:
$$
  \begin{CD}
       X @.\times @. X       @>s>> \!\!    X        \\
   @Af^{-1}AA *@. @AAf^{-1}A       \!\!   @VVfV      \\
       X'@.\times @. X'      @>s'>> \    X'\ .
  \end{CD}
$$

Кроме того, если $f$~--- биекция, то, независимо от того, является
$s$ функцией или нет:
\begin{align*}
  ((x',\ y'),\ z') \in s' &\eq  (\exists x, y, z))(x' = f(x) \land
                                        y' = f(y) \land
                                        z' = f(z) \land
                                       ((x,\ y),\ z) \in s)    \\
  &\eq (\exists x, y, z)(x = f^{-1}(x') \land
                                        y = f^{-1}(y') \land
                                        z = f^{-1}(z) \land
                                       ((x,\ y),\ z) \in s)    \\
  &\eq (f^{-1}(x'),\ f^{-1}(y'),\ f^{-1}(z')) \in s.
\end{align*}

\smallskip
7. \ul{$n = 2,\ T = pr^2_1.$} Для любых $(X,\ Y)$, $(X',\ Y') \in \Scal^2$,
$f: X\to X'$ и $g: Y\to Y'$ имеем:
$$\begin{array}{l}
   T(X,\ Y) = X,\ T(X',\ Y') = X',\\
   T\< f,\ g\> = f.
  \end{array}
$$
Каждой $s \in T(X,\ Y)$ отображение $T\<f,\ g\>$ ставит в соответствие
значение $s'$ функции $f$ в точке $s$.

\smallskip
8. \ul{$n = 2,\ T = pr^2_1\otimes pr^2_2.$} Для любых $(X,\ Y)$, $(X',\ Y')
\in \Scal^2$, $f: X\to X'$ и $g: Y\to Y'$ имеем:
$$\begin{array}{l}
   T(X,\ Y) = X\times Y,\ T(X',\ Y') = X'\times Y',\\
   T\< f,\ g\> = f*g.
  \end{array}
$$
Каждой $s = (x,\ y) \in T(X,\ Y) = X\times Y$ отображение $T\<f,\ g\>$
ставит в соответствие $s' = (f(x),\ g(y)) \in (X',\ Y').$

\smallskip
9. \ul{$n = 1,\ T = \Pcal\circ(pr^2_1\otimes pr^2_1)$.} Этот пример
аналогичен примеру 5. Чтобы понять это, нужно только заметить,
что тождественное отображение $\Dl: \Scal\to \Scal$ есть то же самое,
что и $pr^1_1.$ Для любых
$$\begin{array}{l}
   (X,\ Y),\ (X',\ Y') \in \Scal^2,\\
   f: X\to X',\ g: Y\to Y'
  \end{array}
$$
имеем:
\begin{align*}
   T(X,\ Y)   &= \Pcal(X\times Y),\\
   T(X',\ Y') &= \Pcal(X'\times Y'),\\
   T\<f,\ g\> &= \Pcal\<f * g\>.
\end{align*}
Каждой $s \in T(X,\ Y) = \Pcal(X\times Y)$ отображение $T\<f,\ g\>$
ставит в соответствие
\begin{align*}
   s' &= T\<f,\ g\>(s) = (f*g)[s]\\
      &= \{(f(x),\ g(y)): (x,\ y) \in s\}.
\end{align*}

Отсюда видно, что
\begin{align*}
   (x',\ y') \in s' &\eq (\exists x,\ y)(x' = f(x)     \land y' = g(y) \land (x,\ y) \in s) \\
                    &\eq (\exists x,\ y)((x,\ x')\in f  \land (y,\ y')\in g
                                      \land (x,\ y)\in s) \\
                    &\eq (\exists x,\ y)((x',\ x)\in f^{-1} \land (x,\ y)\in s
                                          \land (y,\ y') \in g) \\
                    &\eq (x',\ y') \in g\circ s\circ f^{-1}.
\end{align*}
Таким образом, $s' = g\circ s\circ f^{-1}.$ Эта связь между $s$ и $s'$
становится особенно прозрачной, если $s$ и $f^{-1}$~--- функции (т. е.
$s$~--- функция и $f$~--- обратимая функция). Тогда $f^{-1}$ и $s'$~---
тоже функции, и соотношение между $s'$, $f^{-1}$, $s$ и $g$ может
быть представлено в виде следующей коммутативной диаграммы:
$$
  \begin{CD}
       X       @>s>>\!\!\! Y        \\
   @Af^{-1}AA       \!\!\! @VVgV      \\
       X'      @>s'>> \  Y'\ .
  \end{CD}
$$
(Коммутативность этой диаграммы состоит в том, что
$s' = g\circ s\circ f^{-1}.$)

Кроме того, если $f$ и $g$~--- биекции, то, независимо от того, является
$s$ функцией или нет:
$$\begin{array}{rcl}
   (x',\ y') \in s' &\eq& (\exists x,\ y)(x' = f(x) \land
                                        y' = g(y) \land
                                       (x,\ y) \in s)    \\
                   &\eq& (\exists x,\ y)(x = f^{-1}(x') \land
                                        y = g^{-1}(y') \land
                                       (x,\ y) \in s)    \\
                   &\eq& (f^{-1}(x'),\ g^{-1}(y')) \in s.
   \end{array}
$$

Если $X=Y=X'=Y'=\R,$ $f = \lm_x(x+1)$ и $g = \lm_x(x - 2)$,
то множества $s'$ получаются из множеств $s\in T(X,\ Y)$ сдвигом
на одну единицу вправо и на две единицы вниз.

\smallskip
10. \ul{$n = 1,\ T = \Pcal\circ((pr^3_1\otimes pr^3_2)\otimes pr^3_3)$.}
Этот пример аналогичен примеру 6. Для любых
$$\begin{array}{l}
   (X,\ Y,\ Z),\ (X',\ Y',\ Z,) \in \Scal^3,\\
   f: X\to X',\ g: Y\to Y',\ h: X\to Z',
  \end{array}
$$
имеем:
$$\begin{array}{l}
   T(X,\ Y,\ Z) = \Pcal((X\times Y)\times Z),\\
   T(X',\ Y',\ Z') = \Pcal((X'\times Y')\times Z'),\\
   T\<f,\ g,\ h\> = \Pcal\<(f * g) * h\>.
  \end{array}
$$
Так же, как в предыдущем примере, убеждаемся, что каждой
$s \in T(X,\ Y,\ Z)$ отображение $T\<f,\ g,\ h\>$ ставит в соответствие
$$
  s' =  \{(f(x),\ g(y),\ h(z)): (x,\ y,\ z) \in s\},
$$
так что $s' = h\circ s\circ (f^{-1}*g^{-1}).$
Снова, если $s$, $f^{-1}$ и $g^{-1}$~--- функции
(т. е. $s$~--- функция и $f$, $g$~--- обратимые функции), то эта
связь между $s$ и $s'$ получает наглядное отражение в коммутативной
диаграмме, аналогичной тем, что приведены выше:
$$
  \begin{CD}
       X @.\times @. X       @>s>> \!\!    X        \\
   @Af^{-1}AA *@. @AAg^{-1}A       \!\!   @VVhV      \\
       X'@.\times @. X'      @>s'>> \    X'\ .
  \end{CD}
$$
А если $f$, $g$ и $h$~--- биекции, то, независимо от того, является
$s$ функцией или нет:
$$
  ((x',\ y'),\ z') \in s' \eq (f^{-1}(x'),\ g^{-1}(y'),\ h^{-1}(z')) \in s.
$$

\section{\bf Перенос структур вдоль биекций}

В предыдущем параграфе мы для каждой последовательности отображений
$f_i\colon X_i\to X'_i,$ $i = 1,\ \dots,\ n,$ и каждого типа ступени $T$ над
$n$ множествами определили отображение $T\<f_1,\ \dots,\ f_n\>$
множества $T(X_1,\ \dots,\ X_n)$ всех структур типа $T$ на множествах
$X_1,\ \dots,\ X_n$ во множество $T(X_1',\ \dots,\ X_n')$ всех структур
типа $T$ на множествах $X_1',\ \dots,\ X_n'.$ Таким образом, задание
отображений $f_i$ позволяет каждой структуре $s\in T(X_1,\ \dots,\ X_n)$
поставить в соответствие структуру $s'\in T(X_1',\ \dots,\ X_n')$,
образно говоря~--- перенести структуру $s$ с множеств $X_1,\ \dots,\ X_n$
на множества $X_1',\ \dots,\ X_n'$.

Осуществляя перенос структур, следует учитывать, что тип структуры
не определяется однозначно. Одно и то же множество может рассматриваться
как структуры различных типов на одних и тех же множествах. Как структуре
одного типа ему при переносе будет отвечать одно множество, как структуре
другого типа~--- вообще говоря, другое множество.

Например, пусть $X = \set{\0}.$ Тогда $\0\in X$, поэтому
$\0$ является структурой типа $pr^1_1$ на $X$. Но, кроме того,
$\0\sbs X$, так что $\0$ является ещё и структурой
типа $\Pcal$ на том же множестве $X$. Пусть теперь $Y = \{\{\0\}\}$
и $f\colon X\to Y$ единственному элементу множества $X$ ставит в соответствие
единственный элемент множества $Y$: $f(\0) = \{\0\}.$
При $T = pr^1_1$ имеем: $T\<f\>(\0) = f(\0) = \{\0\},$ а при $T = \Pcal$~---
$T\<f\>(\0) = \Pcal\<f\>(\0) = f[\0] = \0$. Как видим, результаты
различны.

Итак, совершая перенос структуры, следует указывать, структурой какого
типа она при этом считается. Удобнее всего это делать с помощью формулы
вида `$s\in T(X_1,\ \dots,\ X_n)$', где на месте буквы `$T$', разумеется, стоит
конкретное знакосочетание, изображающее тот или иной тип структуры. Каждую
формулу такого вида будем называть {\it типизацией} буквы `$s$'. Про
структуру $s'$, связанную с $s$ соотношением $s' = T\<f_1,\ \dots,\ f_n\>(s)$,
будем говорить, что она {\it получена из структуры $s$ переносом вдоль
отображений $f_1,\ \dots,\ f_n$ при типизации `$s\in T(X_1,\ \dots,\ X_n)$'}.
Как обычно, если из контекста будет ясно, какая типизация имеется в виду,
указание на неё будет опускаться.

Перенос структур играет исключительно важную роль в математике. Это станет
особенно ясным в следующем параграфе, когда будет определено понятие рода
структуры. Однако уже теперь можно проиллюстрировать роль переноса структур
в математике на хорошо знакомых читателю примерах.

Пусть $s$~--- операция сложения на множестве $\R$ всех действительных чисел,
а $E$~--- экспонента с основанием $e$, то есть биекция множества $\R$
на множество $\R^*_+$, действующая по правилу $E(x) = e^x$, $x\in \R.$
Тогда $s$~--- структура типа $\Pcal\circ((\Dl\otimes\Dl)\otimes\Dl)$ на
множестве $\R$ и как таковой ей при переносе вдоль биекции $E$ отвечает
операция $s' = T\<E\>(s).$ Как было установлено выше (часть 3, пример 6),
$s' = E\circ s\circ (E^{-1}*E^{-1}).$ Поэтому для любых $x,\ y \in \R^*_+$
имеем:
$$
  s'(x,\ y) = E\circ s\circ (E^{-1}*E^{-1})(x,\ y)
           = e^{\ln x + \ln y} = x\cdot y.
$$
Это означает, что операция умножения положительных действительных чисел
получается из операции сложения действительных чисел переносом вдоль
экспоненты (при указанной типизации). Именно на этом соотношении между
указанными операциями основана польза от применения логарифмов. Таким
образом, в основе применения логарифмов для нужд вычислений лежит
идея переноса структур.

Другой пример весьма интенсивного использования переноса структур
доставляет координатный метод. Действительно, система координат
на множестве $X$~--- это просто биекция $x\colon X\to U$ множества $X$
на некоторое множество $U$ подходящего пространства, скажем, $\R^n.$
(Композиции $pr^n_i\circ x$ обычно обозначают через $x_i$ или $x^i$
и называют координатными функциями системы координат $x$.) Задание
системы координат $x$ порождает серию отображений $T\<x\>\colon T(X)\to
T(U)$~--- по одному для каждого типа ступени $T$. Это позволяет
с каждым <<геометрическим образом>>\ (то есть, попросту говоря,
с каждой структурой) $s$ на множестве $X$ соотнести <<такую же>>\
структуру $s' = T\<x\>(s)$ на множестве $U$ и изучение $s$ свести
к изучению $s'$.

Скажем, изучение преобразований $t$ множества $X$ заменяется изучением
преобразований $t'$ множества $U$, связанных с $t$ равенствами вида
$t' = T\<x\>(t).$ Но в данном случае $t' = x\circ t\circ x^{-1}$
(см. часть 3, пример 5), и мы видим, что $t'$ есть то, что обычно
называют координатной формой преобразования~$t$.

Переход от одной системы координат к другой также сопряжён с переносом
структур. В самом деле, всякая замена системы координат на множестве $X$
сводится к выбору системы координат $u\colon U\to U_1$ на множестве $U$. Это
приводит к новой системе координат $u\circ x\colon X\to U_1$ на множестве $X$.
В ней каждой структуре $s$ типа $T$ на множестве $X$ отвечает структура
$s' = T\<u\circ x\>(s)$ на множестве $U_1.$ Но $T\<u\circ x\> =
T\<u\>\circ T\<x\>$ (Т3.3, п. $2^\circ$). Поэтому $s' = T\<u\circ x\>(s) =
(T\<u\>\circ T\<x\>)(s) = T\<u\>(T\<x\>(s)).$ В частности, для
преобразований получаем:
$$
  t' = (u\circ x)\circ t\circ (u\circ x)^{-1}
     = u\circ(x\circ t\circ x^{-1})\circ u^{-1}, \mbox{ ---}
$$
формула, выводу которой в некоторых руководствах посвящаются целые
страницы.

Разумеется, для того, чтобы свойства, обнаруженные у $s',$ можно было
приписать~$s$ (или обратно), нужно, чтобы эти свойства при переносах
<<сохранялись>>. Примеры, рассмотренные в предыдущем параграфе,
показывают, что так бывает не всегда: не всегда топология переходит
в топологию, отношение порядка~--- в отношение порядка и~т.~д.
Однако, оказывается, что если ограничиться рассмотрением переносов
лишь вдоль {\it биекций}, то многие свойства структур уже будут
сохраняться. Свойства структур, сохраняющиеся при переносах вдоль
биекций, называют переносимыми. В следующем параграфе понятию
переносимости будет дано точное определение, а пока мы проиллюстрируем
его на примерах наиболее употребительных типов ступеней и наиболее
часто встречающихся свойств структур.

\begin{thm}[Теорема о переносах бинарных отношений]\sl
Пусть $s$~--- бинарное отношение на множестве $X$, а $s'$ получается из $s$
переносом вдоль биекции $f\colon X\to X'$. Тогда:
\begin{flalign*}
&0^\circ  & (x_1,\ x_2)\in s            &\eq (f(x_1),\ f(x_2))\in s', & \\
&         & (x_1',\ x_2')\in s'         &\eq (f^{-1}(x_1'),\ f^{-1}(x_2'))\in s;\\
&1^\circ  & s' \mbox{ рефлексивно }     &\eq s \mbox{ рефлексивно};\\
&2^\circ  & s' \mbox{ иррефлексивно }   &\eq s \mbox{ иррефлексивно};\\
&3^\circ  & s' \mbox{ симметрично }     &\eq s \mbox{ симметрично};\\
&4^\circ  & s' \mbox{ асимметрично }    &\eq s \mbox{ асимметрично};\\
&5^\circ  & s' \mbox{ антисимметрично } &\eq s \mbox{ антисимметрично};\\
&6^\circ  & s' \mbox{ транзитивно }     &\eq s \mbox{ транзитивно, и т. д.\ \ }
\end{flalign*}
\end{thm}

\begin{proof} П.~$0^\circ$ доказан в части 3, остальные доказываются
простым применением п.~$0^\circ$. Докажем, например,~п.~$5^\circ$. Пусть
$s$\/ антисимметрично, т.~е.
$$
  (\forall x,\ y \in X)((x,\ y) \in s \land (y,\ x) \in s \to x = y).
$$
Покажем, что и $s'$ антисимметрично, т.~е.
$$
  (\forall x',\ y' \in X')((x',\ y') \in s' \land (y',\ x') \in s' \to x' = y').
$$
Пусть
$$
  (x',\ y') \in s' \land (y',\ x') \in s'.
$$
Тогда, в силу п. $0^\circ$,
$$
  (f^{-1}(x'),\ f^{-1}(y'))\in s \mbox{ и } (f^{-1}(y'),\ f^{-1}(x'))\in s.
$$
Теперь антисимметричность $s$ даёт, что $(f^{-1}(x') = f^{-1}(y'))$,
а биективность $f$~--- что $x' = y'.$ Это и требовалось. Обратно~--- точно
так же.
\end{proof}

\begin{thm}[Теорема о переносах бинарных операций]\sl
Пусть $s$ и $t$~--- бинарные операции на множестве $X$, а $s'$ и $t'$ получаются
из $s$ и $t$ переносом вдоль биекции $f\colon X\to X'$. Тогда:
\begin{align*}
& \hspace{-3em} 0^\circ & ((x_1,\ x_2),\ x_3) \in s        &\eq ((f(x_1),\ f(x_2)),\ f(x_3))\in s', &\\
& \hspace{-3em}         & ((x_1',\ x_2'),\ x_3')\in s'     &\eq ((f^{-1}(x_1'),\ f^{-1}(x_2')),\ f^{-1}(x_3'))\in s;\\
& \hspace{-3em} 1^\circ & \mbox{s~--- бинарная операция}   &\eq s'\mbox{ --- бинарная операция};
\smallskip\\
& \hspace{-3em} 2^\circ & s \mbox{ ассоциативна }          &\eq s' \mbox{ ассоциативна};
\smallskip\\
& \hspace{-3em} 3^\circ & s \mbox{ коммутативна }          &\eq s' \mbox{ коммутативна};
\medskip\\
& \hspace{-3em} 4^\circ & \left.\parbox{5.8cm}%
             {существует нейтральный
             элемент относительно $s$}\;\right]  &\eq \left[\;\parbox{5.8cm}{существует нейтральный элемент относительно $s'$;}\right.
\bigskip\\
& \hspace{-3em} 5^\circ & \left.\parbox{5.8cm}%
             {для каждого элемента из $X$
              существует обратный
              элемент относительно $s$}\;\right] &\eq \left[\;\parbox{5.8cm}{для каждого элемента из $X'$ существует обратный элемент
                                                                             относительно $s';$}\right.
\bigskip\\
& \hspace{-3em} 6^\circ & \parbox{5.8cm}{$t$ дистрибутивна
              относительно $s$}                  &\eq \ \parbox{5.7cm}{\;$t'$~дистрибутивна~относительно~$s'$.}
\end{align*}\qedhere

\end{thm}

\smallskip
\begin{proof} Пп. $0^\circ$ и $1^\circ$ уже были доказаны выше
(часть 3, пример 6). (Это было установлено для любой структуры типа
$\Pcal\circ((\Dl\otimes\Dl)\otimes\Dl)$, так что верно и для $s$, и для~$t$.)
Отсюда легко выводятся остальные утверждения теоремы. Докажем, например,
п. $6^\circ.$ Пусть $s$ дистрибутивна относительно $t$, т. е. для любых
$x$, $y$, $z \in X,$
$$
  s(x,\ t(y,\ z)) = t(s(x,\ y),\ s(x,\ z)).
$$
Покажем, что тогда и $s'$ дистрибутивна относительно $t'$, т.~е. для любых
$x'$, $y'$, $z' \in X'$,
$$
  s'(x',\ t'(y',\ z')) = t'(s'(x',\ y'),\ s'(x',\ z')).
$$
При рассмотрении указанного примера (см. приведённую там диаграмму) было установлено, что
$$
  s' = f\circ s\circ(f^{-1}*f^{-1}).
$$
Совершенно аналогично,
$$
  t' = f\circ t\circ(f^{-1}*f^{-1}).
$$
Используя это, дистрибутивность $t$ относительно $s$ и вытекающее из первого равенства (в силу обратимости~$f$) соотношение
$$
  s\circ(f^{-1}*f^{-1}) = f^{-1}\circ s',
$$
для любых $x'$, $y'$, $z' \in X'$ получаем:
\begin{align*}
s'(x',\ t'(y',\ z')) &= (f\circ s\circ(f^{-1}*f^{-1}))(x',\ (f\circ t\circ(f^{-1}*f^{-1}))(y',\ z')) \\
                     &= (f\circ s)(f^{-1}(x'),\ f^{-1}((f\circ t\circ(f^{-1}*f^{-1}))(y',\ z')))\\
                     &= (f\circ s)(f^{-1}(x'),\ t(f^{-1}(y'),\ f^{-1})(z'))\\
                     &= f(s(f^{-1}(x'),\ t(f^{-1}(y'),\ f^{-1})(z')))\\
                     &= f(t(s(f^{-1}(x'),\ f^{-1}(y')),\ s(f^{-1}(x'),\ f^{-1}(z'))))\\
                     &= (f\circ t)((s\circ(f^{-1}*f^{-1}))(x',\ y'),\ (s\circ(f^{-1}*f^{-1}))(x',\ z'))\\
                     &= (f\circ t)((f^{-1}\circ s')(x',\ y'),\ (f^{-1}\circ s')(x',\ z'))\\
                     &= (f\circ t\circ(f^{-1}*f^{-1}))(s'(x',\ y'),\ s'(x',\ z'))\\
                     &= t'(s'(x',\ y'),\ s'(x',\ z')).
\end{align*}
Это и требовалось.
\end{proof}

\begin{thm}[Теорема о переносах топологий]\sl
Пусть $s\in \Pcal(\Pcal(X))$, а $s'$ получается из $s$ переносом вдоль биекции $f\colon X\to X'$. Тогда:
\begin{align*}
&\hspace{-6em} 0^\circ & \mbox{ для любого } U'\sbs X',\ U'\in s'     &\eq f^{-1}[U']\in s; &\\
&\hspace{-6em} 1^\circ & s \mbox{ --- топология}                      &\eq s'\mbox{ --- топология};\\
&\hspace{-6em} 2^\circ & s \mbox{ отделима}                           &\eq s' \mbox{ отделима};\\
&\hspace{-6em} 3^\circ & s \mbox{ связна}                             &\eq s' \mbox{ связна};\\
&\hspace{-6em} 4^\circ & s \mbox{ компактна}                          &\eq s' \mbox{ компактна, и т. д.}
\end{align*}
\end{thm}

\begin{proof} \fbox{1.} Напомним (часть 3, пример 5), что
в рассматриваемом случае $s' = \Pcal\<\Pcal\<f\>\>(s).$ Поэтому
$$
\begin{array}{rcl}
 s' &=& \Pcal\<f\>[s] = \{\Pcal\<f\>(U)\colon U\in s\}\\
    &=& \{U'\colon(\exists U)(U\in s\land U' = f[U])\}\\
    &=& \{U'\colon (\exists U)(U = f^{-1}[U']\land U\in s)\},
\end{array}
$$
Словами: структуре $s'$ принадлежат те и только те подмножества
$U'$ множества $X'$, прообразы которых относительно биекции $f$
принадлежат структуре $s$. Теперь приступим к доказательству п.~$1^\circ.$

Напомним, что топологией на множестве $X$ называется всякая
совокупность его подмножеств, которой принадлежат $\0,$ $X$,
объединения любых частей этой совокупности и пересечения любых
конечных частей этой совокупности.

Покажем, что если $s$~--- топология, то и $s'$~--- топология.
Пусть $s$~--- топология. Тогда:
\begin{itemize}
\item[а)] $\0\in s'$~--- так как $f^{-1}[\0] = \0 \in s;$
\item[б)] $X'\in s'$~--- так как $f^{-1}[X'] = X \in s;$
\item[в)] если $s_1'\sbs s'$, то $\cup s_1'\in s'$~--- так как
$$
  f^{-1}[\cup s'] = \cup\{f^{-1}[U']\colon U'\in s'\} \in s;
$$
\item[в)] если $U_1',\ \dots,\ U_n'\in s'$, то $U_1'\cap\dots\cap U_n'
\in s'$~--- по аналогичным причинам.
\end{itemize}
Всё это и означает, что
$s'$~--- топология. Точно так же доказывается обратное утверждение.
\end{proof}

Доказательства остальных пунктов теоремы оставляются в качестве
упражнений читателям, понимающим, о чём в них идёт речь.

\smallskip
Мы рассмотрели несколько видов структур, определённых на одном
множестве. В математике рассматриваются и структуры, задаваемые
на нескольких множествах. Чаще всего бывает два множества, одно
из которых называется основным, а другое~--- вспомогательным.
На вспомогательном множестве обычно имеется структура, определение
которой не требует обращения к основному множеству, а структура
на основном множестве оказывается некоторым образом связанной
со структурой на вспомогательном множестве. Типичные примеры~---
структура векторного пространства и структура метрического
пространства.

При переносах структур вспомогательные множества ведут себя
особенным образом. Пусть, например, $s$~--- структура векторного
пространства над полем $k$, заданная на множестве $E$. Тогда
$E$~--- основное множество, $k$~--- вспомогательное, а $s$
есть упорядоченная пара отображений:
$$
  a\colon E\times E \to E\quad \mbox{и}\quad m\colon k\times E \to E,
$$
где первое есть операция сложения векторов, второе~--- операция
умножения их на скаляры. Мы можем, взяв произвольную пару биекций
$f\colon E \to E',\ g\colon k \to k',$ перенести $s$ вдоль $(f,\ g).$
В результате на $E'$ образуется структура векторного пространства
над полем $k'$. При таком переносе множество $k$ ведёт себя
<<на равных>>\ с множеством $E.$ Однако такие переносы
в стандартных рассмотрениях теории векторных пространств практически
не встречаются. Совершая перенос структуры векторного пространства
с одного множества на другое, обычно меняют множество векторов,
но не поле скаляров. Иными словами, структуру векторного пространства
обычно переносят вдоль не произвольных пар биекций, а лишь пар вида
$(f,\ I_k),$ где вспомогательному множеству отвечает тождественное
отображение.

Такая же картина наблюдается в теории метрических пространств.
Перенося метрику с одного множества на другое, обычно меняют множество
точек метрического пространства, но не множество, в котором принимает
значения метрика. Таким образом, и здесь рассматриваются переносы лишь
вдоль пар вида $(f,\ I),$ где вспомогательному множеству отвечает
тождественное отображение.

Пусть $T$~--- тип ступени над $n+m$ множествами, $s$~--- структура
типа $T$ на множествах $X_1$, \ldots, $X_n$, $Y_1$, \ldots, $Y_m$,
$f_1\colon X_1 \to X_1'$, \ldots, $f_n\colon X_n\to X_n'$~--- биекции и
$s' = T\<f_1,\ \dots,\ f_n; I_{Y_1},\ \dots,\ I_{Y_m}\>(s).$ Тогда
говорят, что {\it структура $s'$ получена переносом структуры $s$
вдоль биекций $f_1,\ \dots,\ f_n$ при типизации <<$s\in
T(X_1,\ \dots,\ X_n,\ Y_1,\ \dots,\ Y_m)$>>\ с основными базисными
множествами $X_1,\ \dots,\ X_n$ и вспомогательными базисными
множествами $Y_1,\ \dots,\ Y_m.$} Если из контекста ясно, о какой
типизации и при каких основных и вспомогательных множествах
идёт речь, то указания на это опускают.

Отметим, что не следует поддаваться звучанию слова <<вспомогательный>>\
и думать, будто вспомогательные множества играют какую-нибудь
второстепенную или несущественную роль. Их роль не менее важна, чем
роль множеств, называемых основными, а слово <<вспомогательный>>\
указывает лишь на то, что при переносах структур эти множества и их
<<собственные>>\ структуры не меняются.

Предлагаем читателю сформулировать и доказать для структуры векторного
пространства и структуры метрического пространства теоремы, аналогичные
теоремам Т4.1 -- Т4.3.

\section{\bf Род структуры}

Мы почти готовы дать определение главного понятия теории структур~---
понятия рода структуры. Осталось лишь уточнить играющее весьма важную
роль в этом определении понятие переносимости. В предыдущем параграфе
мы уже сказали, что свойство структуры называется переносимым, если оно
<<сохраняется>>\ при переносах вдоль биекций. Теперь нам предстоит
дать точное определение.

Мы видели, что рефлексивность бинарных отношений сохраняется при
переносах вдоль биекций. Точная формулировка этого результата гласила,
что всякий раз, как $s$ является бинарным отношением на множестве $X$,
$f$~--- биекцией $X$ на $X'$ и $s'$ получается переносом $s$ вдоль
биекции $f$, так $s$ рефлексивно тогда и только тогда, когда $s'$
рефлексивно. Таким образом, переносимость свойства рефлексивности
бинарных отношений состоит в том, что является теоремой теории
множеств формула:
\begin{multline*}
  [(s\in\Pcal(X\times X))
     \land (f\mbox{~--- биекция $X$ на $X'$})
     \land (s' = \Pcal\<f*f\>(s))] \to\\
  \to ((s \mbox{~--- рефлексивное отношение на $X$)}
     \eq
      (s' \mbox{~--- рефлексивное отношение на $X'$})).
\end{multline*}
Переносимость свойства транзитивности бинарных отношений означает,
что является теоремой теории множеств формула
\begin{multline*}
  [(s\in\Pcal(X\times X))
     \land (f\mbox{~--- биекция $X$ на $X'$})
     \land (s' = \Pcal\<f*f\>(s))] \to \\
  \to ((s \mbox{~--- транзититивное отношение на $X$)}
     \eq
      (s' \mbox{~--- транзитивное отношение на $X'$})).
\end{multline*}

И т. д. Вообще, переносимость некоторого свойства бинарных отношений,
то есть структур типа $\Pcal\circ(\Dl\otimes\Dl)$, означает, что является
теоремой теории множеств формула вида
$$
  [(s\in\Pcal(X\times X))
     \land (f\mbox{~--- биекция $X$ на $X'$})
     \land (s' = \Pcal\<f*f\>(s))]
  \to ((\dots X \dots s \dots)
     \eq
      (\dots X'\dots s' \dots)),
$$
где <<$\dots X \dots s \dots$>>\ изображает запись на языке теории
множеств того свойства $s$, о переносимости которого идёт речь,
а <<$\dots X' \dots s' \dots$>>~--- то же относительно $s'$.

Переносимость свойства структуры типа $\Pcal\circ((\Dl\otimes\Dl)\otimes\Dl)$
быть бинарной операцией означает, что является теоремой теории множеств
формула
\begin{multline*}
     [(s\in \Pcal((X\times X)\times X))
     \land (f\mbox{~--- биекция $X$ на $X'$})\land (s' = \Pcal\<(f*f)*f\>(s))]\to \\
  \to ((s \mbox{~--- бинарная операция на $X$)}
     \eq
      (s' \mbox{~--- бинарная операция на $X'$})).
\end{multline*}

Переносимость свойства структуры того же типа быть ассоциативной
бинарной операцией означает, что является теоремой теории множеств
формула
\begin{multline*}
     [(s\in \Pcal((X\times X)\times X)
     \land (f\mbox{~--- биекция $X$ на $X'$})\land (s' = \Pcal\<(f*f)*f\>(s))]\to\\
  \to ((s \mbox{~--- ассоциативная бинарная операция на $X$)}
     \eq
      (s' \mbox{~--- ассоциативная бинарная операция на $X'$})).
\end{multline*}
И так далее. Вообще, переносимость некоторого свойства структур
типа $\Pcal\circ((\Dl\otimes\Dl)\otimes\Dl)$ означает, что является
теоремой теории множеств формула вида
\begin{multline*}
     [(s\in \Pcal((X\times X)\times X)
     \land (f\mbox{~--- биекция $X$ на $X'$})\land (s' = \Pcal\<(f*f)*f\>(s))]\to \\
  \to ((\dots X \dots s \dots)
     \eq
      (\dots X'\dots s' \dots)),
\end{multline*}
где <<$\dots X \dots s \dots$>>\ и <<$\dots X' \dots s' \dots$>>\
понимаются так же, как и выше.

Переносимость свойства структуры типа $\Pcal\circ \Pcal$ быть топологией,
компактной топологией и т. д., означает, что является теоремой
теории множеств формула
\begin{multline*}
     [(s\in \Pcal(\Pcal(X)))
     \land (f\mbox{~--- биекция $X$ на $X'$})
     \land (s' = \Pcal\<\Pcal\<f\>\>(s))]\to \\
  \to ((s \mbox{~--- топология на $X$)}
     \eq
      (s' \mbox{~--- топология на $X'$}));
\end{multline*}
\begin{multline*}
[(s\in \Pcal(\Pcal(X)))
     \land (f\mbox{~--- биекция $X$ на $X'$})
     \land (s' = \Pcal\<\Pcal\<f\>\>(s))]\to     \\
   \to ((s \mbox{~--- компактная топология на $X$)}
     \eq
      (s' \mbox{~--- компактная топология на $X'$})).
\end{multline*}
И так далее. Вообще, переносимость какого-либо свойства структуры
типа $\Pcal\circ \Pcal$ означает, что является теоремой теории множеств
формула вида
\begin{multline*}
     [(s\in \Pcal(\Pcal(X)))
     \land (f\mbox{~--- биекция $X$ на $X'$})
     \land (s' = \Pcal\<\Pcal\<f\>\>(s))]
  \to ((\dots X \dots s \dots)
     \eq
      (\dots X'\ldots s' \ldots)).
\end{multline*}
где <<$\dots X \ldots s \ldots$>>\ и <<$\ldots X' \ldots s' \ldots$>>\
понимаются так же, как и выше.

Рассмотренные нами примеры относятся к структурам всего лишь трёх типов
над одним множеством. Однако переход к общему случаю теперь уже не вызывает
затруднения. Довольно ясно, что переносимость какого-либо свойства структуры
типа $T$ над $n$ множествами означает, что является теоремой теории множеств
формула вида
\begin{multline*}
[(s\in T(X_1,\ \ldots,\ X_n))\land (f_1\mbox{~--- биекция $X_1$ на $X_1'$})\land\ldots
\land (f_n\mbox{~--- биекция $X_n$ на $X_n'$})\land (s' = T\<f_1,\ \ldots,\ f_n\>(s))]\to\\
  \to ((\ldots X_1\dots X_n \ldots s \dots)
     \eq
      (\ldots X_1'\ldots X_n' \ldots s' \ldots)).
\end{multline*}
Это~--- если нет вспомогательных множеств. Если же они есть, то
появляются незначительные отличия: в качестве биекций, соответствующих
вспомогательным множествам, берутся тождественные отображения,
причём они в перечне биекций не упоминаются. Дадим теперь точную
формулировку.

Пусть $\a_1$, \ldots, $\a_n$, $\dl$ и $\a_1'$, \ldots, $\a_n'$, $\dl'$ $(n\ges 1)$~---
попарно различные переменные, $\mu_1$, \ldots, $\mu_n$~--- термы,
в которых не встречаются эти переменные, и $\tau$~--- терм, изображающий
тип ступени над $n+m$ множествами. Положим
$$
  \chi = \lc \dl\in\tau(\a_1,\ \ldots,\ \a_n;\ \mu_1,\ \ldots,\ \mu_m)\rc
$$
и назовём $\chi$ типизацией переменной $\dl.$
Пусть, наконец, $\gamma_1,$ \ldots, $\gamma_n$~--- переменные,
отличные  от $\a_1$, \ldots, $\a_n$, $\dl$ и $\a_1'$, \ldots,\ $\a_n'$, $\dl'$
и всех переменных, встречающихся в термах $\mu_1$, \ldots, $\mu_n$, $\tau$.
Говорят, что {\it формула $\vp$ переносима при типизации $\chi$, в которой
буквы $\a_1$, \ldots, $\a_n$ изображают основные базисные множества,
а термы  $\mu_1$, \ldots, $\mu_n$~--- вспомогательные базисные множества},
если является теоремой теории множеств формула
$$\left\{
   \begin{aligned}
    &\dl\in\tau(\a_1, \ldots, \a_n;\ \mu_1, \ldots, \mu_m)\land \\
    &\land \gamma_1 \mbox{~--- биекция $\a_1$ на $\a_1'$}\land\ \dots\\
    &\land \gamma_n \mbox{~--- биекция $\a_n$ на $\a_n'$}\land \\
    &\land \dl' = \tau\<\a_1, \ldots, \a_n;\ I_{\mu_1}, \ldots, I_{\mu_n}\>(\dl)
   \end{aligned}
\right\}
   \to (\vp \eq \vp').
$$

Переходим к определению рода структуры.

\medskip
{\it Род структуры}~--- это текст $\Sigma,$ включающий в себя:
\begin{enumerate}
\item последовательность $\a_1,\ \dots,\ \a_n$ переменных языка теории
множеств, о которых говорят, что они изображают основные базисные
множества;
\item последовательность $\mu_1,\ \dots,\ \mu_m$ термов теории
множеств, о которых говорят, что они изображают вспомогательные
базисные множества;
\item типизацию $\chi = \lc\dl\in\tau(\a_1, \ldots, \a_n;\ \mu_1,
      \ldots, \mu_m)\rc$, где $\dl$~--- переменная, не встречающаяся\\
      в $\lc\tau(\a_1, \ldots, \a_n;\ \mu_1, \ldots, \mu_m)\rc$;
\item формулу $\vp$, переносимую при типизации $\chi$, в которой
      переменные $\a_1, \dots, \a_n$ изображают основные, а термы
      $\mu_1, \ldots, \mu_n$~--- вспомогательные базисные множества.
\end{enumerate}
Формула $\vp$ называется {\it аксиомой} рода структуры $\Sigma$.
{\it Теорией рода структуры}\/ $\Sigma$ называется теория
$\Tcal_{\Sigma}$, получаемая присоединением формулы
$\lc\chi\land\vp\rc$ к аксиомам теории множеств.

\renewcommand{\refname}{\bf Библиографический список}

\end{document}